\documentclass[12pt]{article}
\begin{document}
\begin{center}

{\large\textbf{Preopen sets in bispaces}}\\

\vspace{0.15in}{\textbf{ Amar Kumar Banerjee$^{1}$  and}}
{\textbf{Pratap Kumar Saha$^{2}$}}\\

\vspace{0.15in}\textbf{Abstract}\\
\end{center}

\indent The notion of preopen sets and precontinuity in a topological space was introduced by Mashhour et. al in 1982 [13]. Later the same was studied in a bitopological space in [7] and [9]. Here we have studied the idea of pairwise preopen sets (semi preopen) and pairwise precontinuity (semi precontinuity) in a more general structure of a bispace and investigate how far several results as valid in a bitopological space are affected in a bispace. \\

\vspace{0.05in}\noindent\textbf{2010 AMS Mathematics Subject Classifications :} 54A05, 54E55, 54E99.\\


\noindent\textbf{1.Introduction}\\
\indent The notion of topological space was generalized to a bitopological space by J.C.Kelly [8] in 1963. Levine [12] introduced the idea of semi open sets and semi continuity and Mashhour et al. [13] introduced the concept of preopen sets and precontinuity in a  topological space.  Many works of generalization on bitopological spaces have been seen in [5], [14], [15] etc. Jelic [6] generalized the idea of preopen sets and precontinuity in bitopological space. Later Khedr et al.[9] and A.Kar et al. [7]  further studied the same  in a bitopological space.\\
The idea of a topological space was generalized to a $\sigma$-space (or simply space) by A.D. Alexandroff [1] in 1940 weakening the union requirements where only contable union of open sets were taken to be open. Later Lahiri and Das [11] gave the idea of a bispace generalizing the notion of bitopological spaces. The concept of semi open sets and quasi open sets in the setting of more general structure of a bispace were studied in [2], [3]. In this paper we wish to study the idea of preopen sets, precontinuity, semi preopen sets and semi precontinuity in  more general structure of a bispace. We have also investigated here how far several results as valid in bitopological space are affected in a bispace.\\

\noindent\textbf{2. Preliminaries :}\\
 \noindent\textbf{Definition 2.1[1]:} A set $X$ is called an Alexandroff space or simply a space if in it is chosen a system $\mathcal{F}$ of subsets  satisfying the following axioms:\\ \indent 1. The intersection of a countable number of sets from $\mathcal{F}$ is a set in $\mathcal{F}$.\\ \indent 2. The union of a finite number of sets from $\mathcal{F}$ is a set in $\mathcal{F}$.\\ \indent 3. The void set $\phi$ is a set in $\mathcal{F}$.\\ \indent 4. The whole set $X$ is a set in $\mathcal{F}$.\\
\indent Sets of $\mathcal{F}$ are called closed sets. Their complementary sets are called open . It is clear that instead of closed sets in the definition of the space, one may put open sets with subject to the conditions of countable summability, finite intersectibility and the condition that $X$ and void set  $\phi$ should be open. The collection of all such open sets will sometimes be denoted by $\tau$ and the space by $(X,\tau)$. Note that a topological space is a space but in general $\tau$ is not a topology as can be easily seen by taking $X=R$, the set of  real numbers and $\tau$ as the collection of all $F_{\sigma}$ -sets in $R$.\\
\noindent\textbf{Definition 2.2[1]:} To every set $M$ of a space $(X,\tau)$ we correlate its closure $\overline{M},$ the intersection of all closed sets containing $M$. The closure of a set $M$ will be denoted by $\tau clM$ or simply $clM$ when there is no confusion about $\tau$.  \\
\indent Generally the closure of a set in a space may not be a closed set. The definition of limit point of a set is parallel as in the case of a topological space.\\    \indent From the axioms, it easily follows that \\ \indent $1) \overline{M\cup N}=\overline{M}\cup \overline{N}$; \indent $2) M\subset \overline{M}$ ; \indent $3) \overline{M}=\overline{\overline{M}}$ ; \indent $4) \overline{\phi} = \phi$; \indent $5) \overline{A}=A\bigcup A^{'}$ where $A^{'}$ denotes the set of all limit points of $A$.\\
\noindent\textbf{Definition 2.3[10]:} The interior of a set $M$ in a space $(X,\tau )$ is defined as the union of all open sets contained in $M$ and is denoted by $\tau - intM$ or $intM$ when there is no confusion about $\tau$.\\
\noindent\textbf{Definition 2.4[8]}: A set $X$ on which are defined two arbitrary topologies is called a bitopological space and in denoted by $(X,P,Q)$ .\\
\noindent\textbf{Definition 2.5[11]:} Let $X$ be a nonempty set. If $\tau_{1}$ and $\tau_{2}$ be two collections of subsets of $X$ such that $(X,\tau_{1})$ and $(X,\tau_{2})$ are two spaces, then $X$ is called a bispace and is denoted by $(X,\tau_{1}, \tau_{2}) $.\\
\noindent\textbf{Definition 2.6[13]:} Let $(X,\tau )$ be a topological space. A subset $A$ of $X$ is said to be preopen if $A\subset int(cl A).$\\
\noindent\textbf{Definition 2.7[9]:} A subset $A$ of bitopological space $(X,\tau_{1}, \tau_{2}) $ is said to be $(\tau_{i},\tau_{j})$ preopen ( or briefly $(i,j)$  preopen ) if there exists $U\in \tau_{i}$ such that $A\subset U\subset \tau_{j}cl( A )$ or equivalently $A\subset \tau _{i} int (\tau _{j} cl(A)).$\\
\indent A subset $A$ is said to be pairwise preopen if it is $(i,j)$ preopen for $i,j=1,2; i\neq j$.\\

\noindent\textbf{3. Preopen  sets:}\\
\indent Throughout our discussion,  $(X,\tau_{1},\tau_{2})$ or simply $X$ stands for a bispace, R stands for the set of real numbers , Q stands for the set of rational numbers, P for the set of irrational numbers, N for the set of natural  numbers and sets are always subsets of $X$ unless otherwise stated.\\
\indent In a topological space the following conditions (I) and (II) are equivalent.\\
Condition (I) : $A\subset int(cl(A))$ and \\
Condition (II) : there exists an open set $U$ such that $A\subset U \subset cl (A).$\\
\indent But in a space (Alexandroff space ) these two condition are not equivalent. In fact condition (I) is weaker than the condition (II) as shown in the example 3.1. In view of above we take condition (II) to define preopen sets in a space. \\
\textbf{Definition 3.1} : Let $(X,\tau )$ be a space. A subset $A$ of $X$ is said to be preopen if there exists an open set $U$ such that $A\subset U \subset cl (A)$ and $A$ is said to be weakly preopen if $A\subset int(cl(A))$.\\
\indent Note that if $A$ is open then $int A=A$. So $A=intA\subset int(clA)$  and hence $A$ is weakly preopen. Also if $A$ is open then the condition (II) holds if we take $U=A$. \indent So every open set is preopen also. But converse may not be true as shown in the following example.\\
\textbf{Example 3.1} : Let $X=[1,2]$ and $\tau = \{X,\phi , F_{i}\}$ where $F_{i}$'s are the countable subsets of irrational in $[1,2]$. Let $A=[1,2]-Q$. Then $cl(A)=X$. So $A$ is weakly preopen, since condition (I) holds. Also condition (II) holds if we take $U=X$. Hence $A$ is preopen also. But $A$ is not open. Next, let us consider $B=([1,2]-Q)-\{\sqrt{2}\}$. Then $cl(B)=X-\{\sqrt{2}\}$ and $int(cl(B))=B.$ Therefore, $B\subset int(cl(B))$. Again, there does not exists any $\tau$ open set $U$ such that $B\subset U\subset cl(B)$. So $B$ is weakly preopen but not preopen.\\
\indent In the definition of preopen sets in a bitopological space $(X,\tau_{1}, \tau_{2}) $  Khedr et al.[9] used the following condition $C_{1}$ where as Kar et al. [7] used the condition $C_{2}$.  \\
\indent $C_{1}$ : there exists a $U\in \tau_{i}$ such that $A\subset U \subset \tau _{j}cl (A).$\\
\indent $C_{2}$ : $A\subset int_{i}(cl_{j}A)$\\
It can be checked that the conditions $C_{1}$ and $C_{2}$ are equivalent in  a  bitopological  space. But in case of a bispace $(X,\tau_{1}, \tau_{2}) $ it can be easily checked that the condition $C_{1}$ implies the condition $C_{2}$. But converse may not be true as shown in the following example.\\
\textbf{Example 3.2} : Let $X=[0,2]$ and $\tau_{1}=\{X,\phi, F_{i}\}$, $\tau_{2}=\{X,\phi, G_{i}\}$ where $F_{i}$'s and $G_{i}$'s are the countable subsets of $[0,1]-Q$ and $[1,2]-Q$ respectively. Then clearly $(X,\tau_{1}, \tau_{2}) $ is a bispace which is not a bitopological space. Let $A=[0,1]-Q$ then $\tau_{2}clA=[0,1]\cup ([1,2]\cap Q) \neq X$. Since $\tau_{1}$ open sets other than $X$ are countable and $A$ contain uncountable number of irrational numbers, there does not exists any $\tau_{1}$ open set $U$ such that $A\subset U \subset \tau _{2}cl (A)$ holds. But $\tau_{1}int(\tau_{2}clA)=\tau_{1}int([0,1]\cup ([1,2]\cap Q))=[0,1]-Q=A$. Therefore, $A\subset \tau_{1}int(\tau_{2}clA)$.\\
\indent In view of above discussion, we give the definition of pairwise preopenness in a bispace as follows:\\
\textbf{Definition 3.2(cf.[9])} : Let $(X,\tau_{1}, \tau_{2}) $ be a  bispace  and $A$ be a subset of $X$  then $A$  is said to be $(\tau_{i},\tau_{j})$ preopen ( or briefly $(i,j)$  preopen ) if there exists $U\in \tau_{i}$ such that $A\subset U\subset \tau_{j}cl A$ , $i,j=1,2; i\neq j$. \\
\indent $A$  is called pairwise preopen if $A$ is both $(\tau_{1},\tau_{2})$ and $(\tau_{2},\tau_{1})$ preopen.\\
\textbf{Note 3.1} : In [7], it is shown that pairwise preopenness in a bitopological space $(X,\tau_{1}, \tau_{2}) $ does not imply the preopenness in the individual topological spaces $(X,\tau_{1}) $ and $(X, \tau_{2}) $. The same assertion is valid for pairwise preopen sets in a bispace as shown in the following example.\\
 \textbf{Example 3.3} : Let $X=[1,3]$, $\tau_{1}=\{X,\phi, G_{i}\cup \{\frac{5}{2}\}\}$  and $\tau_{2}=\{X,\phi, F_{i}\cup \{\frac{3}{2}\}\}$ where $G_{i}$'s and $F_{i}$'s are countable subsets of irrational numbers in $[1,\sqrt{3}]$ and $[\sqrt{3},3]$ respectively. Then $(X,\tau_{1}, \tau_{2}) $ is a bispace which is not a bitopological space. Let $A=\{\sqrt{3}\}$, then $\tau_{1}clA=(X-P_{1})-\{\frac{5}{2}\}$ where $P_{1}$ is set of all irrational numbers in $[1,\sqrt{3})$ and $\tau_{2}clA=(X-P_{2})-\{\frac{3}{2}\}$ where $P_{2}$ is set of all irrational numbers in $(\sqrt{3},3]$.  Then $U=\{\sqrt{3},\frac{5}{2}\}$ is a $\tau_{1}$ open set and $V=\{\frac{3}{2},\sqrt{3}\}$ is a $\tau_{2}$ open set containing $A$ and $A\subset U\subset \tau_{2}clA$ and $A\subset V\subset \tau_{1}clA$. But there does not exists any $\tau_{1}$ open set  $U$ such that $A\subset U\subset \tau_{1}clA$ holds and also there does not exists any $\tau_{2}$ open set  $V$ such that $A\subset V\subset \tau_{2}clA$ holds. So  $A$ is pairwise preopen but not individually preopen.   \\
\textbf{Note 3.2} : Clearly every $\tau_{i}$ open set is $(\tau_{i},\tau_{j})$ preopen but converse may not be true which is shown in the following example :\\
\textbf{Example 3.4} : Let $X=[0,3]$,  $\tau_{1}=\{X,\phi, F_{i}\cup \{\sqrt{2}\}\}$  and $\tau_{2}=\{X,\phi, G_{i}\}$ where $F_{i}$'s are the countable subsets of rational number in $[0,1]$ and $G_{i}$'s are countable subsets of irrational number in $[2,3]$. Then $(X,\tau_{1}, \tau_{2}) $ is a bispace which is not a bitopological space. Let $A=\{0,1\}$ which is not $\tau_{1}$ open. Clearly $\tau_{2}clA=[0,2]\cup ([2,3]\cap Q)$ and also the set $U=\{0,1\}\cup \{\sqrt{2}\}\in \tau_{1}$. So $A\subset U\subset \tau_{2}clA$.  Thus $A$ is $(\tau_{1},\tau_{2})$ preopen but not $\tau_{1}$ open.\\
\textbf{Definition 3.3[3]} : Let $(X,\tau_{1}, \tau_{2}) $ be a bispace. Then a subset $A$ of $X$ is  said  to  be $\tau_{i}$ semi open with respect to $\tau_{j}$ if and only if there exists a $\tau_{i}$ open set $O$ such that $O\subset A\subset \tau_{j}clO$ $i,j=1,2; i\neq j$\\
\textbf{Definition 3.4 (cf.[9])} : A subset $A$ of a bispace $(X,\tau_{1}, \tau_{2}) $ is said to be $(\tau_{i},\tau_{j})$ semi preopen ( or briefly $(i,j)$ semi  preopen ) if there exists an $(i,j)$ preopen set $U$ such that $U\subset A\subset \tau_{j}cl U$ , $i,j=1,2; i\neq j$. \\
\indent $A$ is called pairwise semi preopen if $A$ is both $(\tau_{1},\tau_{2})$ and $(\tau_{2},\tau_{1})$ semi preopen. \\

\indent Clearly every $\tau_{i}$ semi open set with respect to $\tau_{j}$ is $(i,j)$ semi preopen and similarly every $(i,j)$ preopen set is $(i,j)$ semi preopen but  converse  may  not  be  true  as  shown  in  the  following  example.\\
\textbf{Example 3.5} : Example of a $(i,j)$ semi preopen set which is not $(i,j)$ preopen and $\tau_{i}$ semi open w.r.to $\tau_{j}$.\\
\indent Let $X,\tau_{1}$ and $\tau_{2}$ be  as in example 3.4. Let $B=\{0,1\}\cup \{\sqrt{3}\}$. Then $B$ is not $(\tau_{1},\tau_{2})$ preopen. But $A=\{0,1\}$ is $(\tau_{1},\tau_{2})$ preopen [ by example 3.4 ] and $A\subset B\subset \tau_{2}cl A$. So $B$ is $(\tau_{1},\tau_{2})$ semi preopen. Again, $B$ does not contain any $\tau_{1}$ open set. So $B$ can not be $\tau_{1}$ semi open w.r. to $\tau_{2}$.\\
\textbf{Theorem 3.1} : Let $A$ be  a subset in a bispace $(X,\tau_{1}, \tau_{2}) $. (a) If there exists an $(i,j)$ preopen set $U$ such that $A\subset U\subset \tau_{j}cl A$ then $A$ is $(i,j)$ preopen.  (b) If there exists an $(i,j)$ semi preopen set $U$ such that $U\subset A\subset \tau_{j}cl U$, then $A$ is $(i,j)$ semi preopen.\\
\textbf{Proof  :} (a) : Since $U$ is $(i,j)$ preopen, there exists a $\tau_{i}$ open set $G$ such that $U\subset G\subset \tau_{j}cl U$. So $A\subset U\subset G\subset \tau_{j}cl U\subset \tau_{j}cl(\tau _{j}cl A)=\tau _{j}cl A$. This implies $A$ is $(i,j)$ preopen.\\
 \indent (b) : Next, since $U$ is $(i,j)$ semi preopen there exists $(i,j)$ preopen set $V$ such that $V\subset U\subset \tau_{j}cl V$. Therefore, $V\subset U\subset A\subset \tau_{j}cl U\subset \tau_{j}cl(\tau _{j}cl V)=\tau _{j}cl V$. This implies $A$ is $(i,j)$ semi preopen.\\
\textbf{Theorem 3.2} : In a bispace $(X,\tau_{1}, \tau_{2}) $. If $A\subset X$ is $(i,j)$ preopen then for every $\tau _{j}$ closed set $G$ containing $A$, $A\subset \tau_{i}intG$, $i,j=1,2; i\neq j$.\\
\textbf{Proof  :} Let $A$ be $(i,j)$ preopen. Then there exists an $U\in \tau_{i}$ such that $A\subset U\subset \tau_{j}cl A$. Let $G$ be any $\tau_{j}$ closed set containing $A$. Then $\tau_{j}clA\subset G$ and so $A\subset U \subset \tau_{i}int(\tau_{j}clA)\subset \tau_{i}intG$.\\
\textbf{Note 3.3 :} The converse  of  the  theorem 3.2 is true in a bitopological space as seen in theorem 3.8  of  [7]. But unlikely in case of  a   bispace   the converse may not be true as shown in the following example.\\
\textbf{Example 3.6} : Taking $X, \tau_{1}, \tau_{2}$ and $A$ as in example 3.2. Then every $\tau_{2}$ closed set $G$ contains $A$ and $\tau_{1}intG=A$. So $A\subset \tau_{1}intG$. But $A$ is not $(\tau_{1},\tau_{2})$ preopen.\\
\textbf{Theorem 3.3} : Countable union of $(i,j)$ preopen [$(i,j)$ semi preopen] sets is $(i,j)$ preopen [$(i,j)$ semi preopen].\\
\textbf{Proof} : Let $\{A_{k}: k=1,2,3.....\}$ be a countable collection of $(i,j)$ preopen sets in the bispace $(X,\tau_{1}, \tau_{2}) $. Let $A=\bigcup\limits_{k=1}^{\infty}A_{k}$. Since each $ A_{k}$ is $(i,j)$  preopen there exists $\tau_{i}$ open set $ U_{k}$ such that $A_{k}\subset U_{k}\subset \tau_{j}cl(A_{k}) $. Since $ A_{k}\subset A $ for each $k$, $\tau_{j}cl(A_{k})\subset \tau_{j}cl(A)$ for each $k$. So $A_{k}\subset U_{k}\subset \tau_{j}cl(A) $ for  each  $k$.  This implies that $\bigcup\limits_{k=1}^{\infty}A_{k}\subset\bigcup\limits_{k=1}^{\infty}U_{k}\subset \tau_{j}cl(A) $ i.e., $A\subset U\subset \tau_{j}cl(A)$ where $U=\bigcup \limits_{k=1}^{\infty}U_{k}$ is $\tau_{i}$ open.\\
\indent The proof for the case that $A$ is $(i,j)$ semi preopen is similar.\\
Arbitrary union of $(i,j)$ preopen set may not be $(i,j)$ preopen as shown in the following examples.\\
\textbf{Example 3.7} : Take $X,\tau_{1}$ and $\tau_{2}$ as in example 3.2. Let $A_{s}=\{s\}$ where $s\in [0,1]-Q$. Then $A_{s}$ is $\tau_{1}$ open and hence $(\tau_{1},\tau_{2})$ preopen. But $\cup \{A_{s}:s\in [0,1]-Q\}=[0,1]-Q$ is not $(\tau_{1},\tau_{2})$ preopen as shown in example 3.2.\\
\textbf{Example 3.8} : Let $X=[0,3]$, $\tau_{1}=\{X, \phi , F_{i}\cup \{\frac{3}{2}\}\}$ where $F_{i}$'s are the countable subset of $[0,1]-Q$ and $\tau_{2}=\{X,\phi,G_{i}\}$  where $G_{i}$'s are the countable subset of $[2,3]-Q$. Obviously $( X,\tau_{1},\tau_{2})$ is a bispace  but  not a bitopological space. Now consider a set $A_{s}=\{s\}$ where $s\in [0,1]-Q$.  Then $A_{s}$ is not $\tau_{1}$ open. Now  $\tau_{2}cl(A_{s})=[0,2]\cup ( [2,3]\cap Q)$. Also  $A_{s}\cup\{\frac{3}{2}\}$  is  a $\tau_{1}$ open  set  and $A_{s}\subset A_{s}\cup\{\frac{3}{2}\}\subset \tau_{2}cl(A_{s})$.  Therefore,  $A_{s}$ is $(\tau_{1},\tau_{2})$  preopen  but  not $\tau_{1}$ open.  But $\bigcup \{A_{s} : s \in [0,1]-Q \} =[0,1]-Q$ is  not  $(\tau_{1},\tau_{2})$  preopen,  since  $\tau_{2}cl([0,1]-Q)= [0,2] \cup ( [2,3] \cap Q) \neq X$  and  the  $\tau_{1}$ open  set  containing  $[0,1]-Q$  is  only  the  set $X$.\\
\textbf{Remark 3.1} : In example 3.3  of  [9] it is shown that intersection of  two  $(i,j)$ preopen [$(i,j)$ semi preopen ] set in a bitopological  space  may not be $(i,j)$  preopen [ $(i,j)$ semi preopen ]. So the same  assertion is  true  in a bispace.\\
\textbf{Lemma 3.1 :} In  a  bispace  $(X,\tau_{1}, \tau_{2}) $,  if $A \subset X$ and $B \in \tau _{j}$ then $\tau _{j} cl (A) \cap B \subset \tau _{j} cl (A \cap B)  $, for $j = 1,2$.\\
\textbf{Proof} : Let $x \in \tau _{j} cl (A) \cap B$ then $x \in B$  and $ x \in \tau _{j} cl (A)$  where $\tau _{j} cl (A) = A \cup A ^{\prime^{j}} $ ,  $A^{\prime^{j}}$ beings  the  set  of  all $\tau _{j}$  limit  points  of $A$.  So  $x \in A \cup A^{\prime^{j}}$.\\
Case I : Now  if $x \in A$  then $x \in A \cap B$,  so  $x \in \tau _{j} cl (A \cap B)$.\\
Case II : If $x \in A^{\prime^{j}}$,  then $x$  is  a $\tau _{j}$ limit  point  of $A$.  Let $V$ be  any $\tau _{j}$ open  set  containing $x$  then $V \cap B$ is  also $\tau _{j}$ open  set  containing  $x$.  Since $x$  is  a $\tau _{j}$ limit  point  of $A$,  $( V \cap B ) \cap ( A - \{ x \} ) \neq \phi$ that  is $ V \cap ( (A\cap B) - \{ x \} ) \neq \phi$.   So $x$ is  a $\tau _{j}$  limit  point  of  $A \cap B$. Hence $x \in \tau _{j} cl (A \cap B)$.  Hence  the  result  is  proved.\\
\textbf{Theorem 3.4 :} In  a  bispace  $( X,\tau_{1},\tau_{2})$,   if $A$ is $( i , j )$ preopen [ $( i , j )$ semi preopen] and $B \in \tau _{1} \cap \tau _{2}$ then $A \cap B $ is $( i , j )$  preopen [$( i , j )$ semi preopen ], $i \neq j , i,j = 1,2$.\\
\textbf{Proof :} Since $A$ is $( i , j )$ preopen  there  exists $\tau _{i}$ open  set $U$ such  that $A \subset U \subset \tau _{j} cl (A)$.  Now  $A \cap B \subset U \cap B \subset \tau _{j} cl (A) \cap B \subset \tau _{j} cl (A \cap B) $,  since $B$ is $\tau _{j}$ open.  Again $U \cap B$ is $\tau _{i}$ open,  since $B$  is $\tau _{i}$ open  also.  Therefore,  $A \cap B $ is $( i , j )$ preopen. Next,  let $A$ be $( i , j )$ semi  preopen  set.  So  there  exists a $( i , j )$ preopen  set $O$ such  that $O \subset A \subset \tau _{j} cl (O)$.  So $O \cap B \subset A \cap B \subset \tau _{j} cl (O) \cap B \subset \tau _{j} cl (O\cap B)$,  since $B$  is  $\tau _{j}$ open.  Since  by  first part of the  theorem, $O \cap B$ is $( i , j )$ preopen,   $A \cap B$ is $( i , j )$ semi  preopen.\\
\textbf{Note 3.4 :} Let $Y \subset X$  and  let $\tau _{i/Y} = \{ U \cap Y : U \in \tau _{i} \};  i=1,2.$ Then  $( Y ,  \tau _{1/Y}, \tau _{2/Y} )$ is  a  bispace  called  subbispace  of $( X,\tau_{1},\tau_{2})$.  As  in  the  case  of  a  topological space,  it  can  be  checked  that  if $A \subset Y \subset X$ then  $\tau _{i/Y} cl (A) =\tau _{i} cl (A) \cap Y$.\\
\textbf{Theorem 3.5 :}  If $A \subset Y \subset X$ in  a  bispace  $( X,\tau_{1},\tau_{2})$ and  if  $A$ is $( \tau_{i},\tau_{j} )$  preopen [ $( \tau_{i},\tau_{j} )$ semi preopen] in $X$ then $A$ is $( \tau_{i/Y},\tau_{j/Y} )$  preopen  [ $( \tau_{i/Y},\tau_{j/Y} )$ semi preopen] in $( Y, \tau_{i/Y},\tau_{j/Y} )$  $i , j = 1, 2, i \neq j$.If, in addition,  $Y \in \tau _{i}$,  then  the  converse  hold.\\
\textbf{Proof :}  Let $A$ be $( \tau_{i},\tau_{j} )$  preopen  in $X$  then  there  exists  an $U \in \tau _{i}$ such  that $A \subset U \subset \tau _{j} cl (A)$. Therefore, $A \subset U \cap Y \subset \tau _{j} cl (A) \cap Y = \tau _{j/Y}  cl ( A )$.  Since $U \cap Y$ is $\tau _{i/Y}$ open, $A$ is $( \tau_{i/Y},\tau_{j/Y} )$ preopen.  Conversely  let $A$ be $( \tau_{i/Y},\tau_{j/Y} )$ preopen.  So  there exists $\tau _{i/Y}$ open  set $U$ such  that $A \subset U \subset \tau _{j/Y} cl (A)$. Since  $Y$  is $\tau _{i}$ open, $U$ is $\tau _{i}$ open. Also $\tau _{j/Y} cl (A) =\tau _{j} cl (A) \cap Y \subset \tau _{j} cl (A)$.   So  $A \subset U \subset \tau _{j} cl (A)$ where $U$ is $\tau _{i}$ open.  So  $A$ is  $( \tau_{i},\tau_{j} )$  preopen. The proof  for the case of semi preopenness  is  similar.\\
\textbf{Definition 3.5 (cf. [9] ) :}  A  subset $A \subset X$ is  said  to  be $( \tau_{i},\tau_{j} )$  or $( i , j )$  preclosed [ resp. $( \tau_{i},\tau_{j} )$ or $( i , j )$ semi preclosed ]  if  its  complement  is $( \tau_{i},\tau_{j} )$  or $( i , j )$ preopen [  resp. $( \tau_{i}, \tau_{j} )$  or $( i , j )$ semi preopen ]. $A$ is  called  pairwise  preclosed [ resp. pairwise  semi  preclosed ] if  its  complement  is  pairwise  preopen [resp. pairwise semi preopen ].\\
\textbf{Definition 3.6 ( cf. [9] ) :} The  intersection  of  all $( i , j )$ preclosed [ resp. $( i , j )$ semi preclosed ] sets  containing $A$ in $( X,\tau_{1},\tau_{2})$ is  called  $( i , j )$ preclosure [ resp. $( i , j )$ semi preclosure ] of  $A$ and  is  denoted  by $( i , j ) pcl (A)$ or $( \tau_{i},\tau_{j} ) pcl (A)$ [ resp. $( i , j ) spcl (A)$ or  $( \tau_{i},\tau_{j} ) spcl (A)$ ], $i, j = 1, 2;  i \neq j $.\\
\textbf{ Theorem 3.6 }  Let $A$ and $B$ be  subsets  of  $X$  in  a  bispace  $( X,\tau_{1},\tau_{2})$ and  let  $x \in X $, then  \\
 (i)  $x \in ( i, j ) pcl (A)$ if  and  only  if $A \cap U \neq \phi$ for  every $( i , j )$ preopen  set $U$ containing $x, i, j = 1, 2;  i \neq j $. \\
 (ii) if $A \subset B$  then $( i , j ) pcl (A) \subset ( i, j ) pcl (B)$.\\
\textbf{Proof :} (i) Let  $U$ be  a  $( i , j )$   preopen  set  containing $x$ such  that $A \cap U = \phi$. This  implies  that $A \subset X - U$, where $X - U$  is $( i , j )$ preclosed.  So $( i , j ) pcl (A) \subset X-U$. Since $x \not \in X - U$, $x \not \in ( i , j ) pcl (A) $.  Therefore, if  $x \in ( i , j ) pcl (A)$    then  every $( i , j )$ preopen  set $U$ intersects $A$. Conversely,  let  every $( i , j )$ preopen  set  containing $x$ intersects  $A$ and  let $x \not \in ( i , j ) pcl (A)$. Then  there  exists  a $( i , j )$ preclosed  set $F$ containing $A$ such  that $x \not \in F$. So $X - F$ is  an $( i , j )$ preopen  set  containing $x$ such  that $( X -F ) \cap A = \phi$ which  is  a  contradiction.  Hence $x \in ( i , j ) pcl ( A )$. \\
The  proof of  (ii) is  straight  forward  and so is  omitted.\\
\textbf{ Theorem 3.7 : }  Let $A$ and $B$ be  subsets  of  $X$  in  a  bispace $( X,\tau_{1},\tau_{2})$ and let  $x \in X $, then  \\
\indent (i)  $x \in ( i, j )  spcl (A)$ if  and  only  if $A \cap U \neq \phi$ for  every $( i , j )$ semi preopen  set $U$ containing $x, i, j = 1, 2;  i \neq j $. \\
\indent (ii) if $A \subset B$  then $( i , j ) spcl (A) \subset ( i, j ) spcl (B)$.\\
\textbf{Proof :} Proof  is similar  as  the  proof  of    theorem 3.6.\\

\noindent\textbf{4. Pairwise  precontinuity :} \\
\textbf{Definition 4.1 [2] :} A function  $f : ( X,\tau _{1},\tau _{2} ) \rightarrow ( Y, \sigma _{1}, \sigma _{2} ) $ is said to be  pairwise  continuous  or  simply continuous  if  the  induced  functions $ f : ( X,\tau _{1} ) \rightarrow ( Y, \sigma _{1} ) $  and $ f : ( X,\tau _{2} ) \rightarrow ( Y, \sigma _{2} ) $ both  are  continuous.\\
\textbf{Definition 4.2 ( cf. [9] )  :} A function  $f : ( X,\tau _{1},\tau _{2} ) \rightarrow ( Y, \sigma _{1}, \sigma _{2} ) $ is said to be  pairwise  open  or  simply open  if  the  induced  functions $ f : ( X,\tau _{1} ) \rightarrow ( Y, \sigma _{1} ) $  and $ f : ( X,\tau _{2} ) \rightarrow ( Y, \sigma _{2} ) $ both  are  open.\\
\textbf{Definition 4.3 ( cf. [9] )  :} A function  $f : ( X,\tau _{1},\tau _{2} ) \rightarrow ( Y, \sigma _{1}, \sigma _{2} ) $ is said to be pairwise  precontinuous  if  the inverse  image  of  each  $\sigma _{i}$  open set of $Y$ is $( i , j )$ preopen in $X$ where $i \neq j ; i, j =1,2.$\\
\indent Clearly  every  pairwise  continuous  function is  pairwise  precontinuous. But  converse  may  not  be  true  even  in  the  case  of  a  bitopological space  as  seen  in  example   4.2[9].\\
\textbf{Note 4.1 :}  If $( X, \tau )$ and $( Y, \sigma )$ be  two  $\sigma$-spaces  and $f : X \rightarrow Y$  is a  continuous function  then  as  in  the  case  of  a  topological  space  the    condition   $f ( \tau cl (A) ) \subset \sigma cl ( f (A) )$ holds good. But  converse of  this  result  may  not  hold as shown  in  the  following  example  although  the  converse  also  holds  in  case  of  topological  spaces.\\
\textbf{Example 4.1 :} Let $X = [0,1]$  and $Y = [1,2]$ and let $\tau = \{ \phi, X, F _{i} \}$ and $\sigma = \{ \phi, X, G _{i} \}$ where $F _{i}$'s and $G _{i}$'s are  countable  subsets  of  irrational  numbers  in $X$ and $Y$ respectively. Then $( X, \tau )$  and $( Y, \sigma )$ are two  spaces  but  not  topological  spaces.   Now  consider  a  function $f : ( X, \tau ) \rightarrow ( Y, \sigma )$ such that\\
$f (x) = \sqrt{2}$  if $x$ is  irrational\\
\indent\indent $=\frac{3}{2}$ if  $x$ is  rational,\\
Now  $\{ \sqrt{2} \} \in \sigma$ and $f ^{-1} ( \{ \sqrt{2} \} ) = [ 0, 1]-Q$ which  is  not $\tau$ open. So $f$ is  not  continuous. Let $A$ be  any  subset  of $X$. Then  the  following  cases  may arise.\\
Case I : $A$ contains  only  rational  numbers.  Then $\tau cl (A) = [0,1] \cap Q$ and $f(\tau cl (A)) = \{\frac{3}{2}\}$. Also $f ( A ) = \{\frac{3}{2}\}$ and $\sigma cl (f(A)) = [1,2] \cap Q$. So  $f ( \tau cl (A) ) \subset \sigma cl ( f (A) )$.\\
Case II : $A$ contains  only  irrational  numbers.  Then $\tau cl (A) = ([0,1] \cap Q) \cup A$ and $f(\tau cl (A)) = \{ \frac{3}{2}, \sqrt{2} \}$,  $f ( A ) = \{ \sqrt{2} \}$ and $\sigma cl (f(A)) = ([1,2] \cap Q) \cup \{ \sqrt{2} \}$. So  $f ( \tau cl (A) ) \subset \sigma cl ( f (A) )$.\\
Case III :  $A$ contains  rational  and  irrational  numbers of $X$. Then  let $A = A _{1} \cup A _{2}$ where $A _{1}$ contains  only  rational  points  of $A$ and $A _{2}$ contains  only  irrational  points  of $A$.  So $\tau cl (A) = \tau cl ( A _{1} \cup A _{2} )= \tau cl (A _{1}) \cup \tau cl (A _{2}) $. Now $\tau cl ( A _{1} ) = [0,1] \cap Q$ and $\tau cl ( A _{2} ) = ([0,1] \cap Q) \cup A _{2}$. Therefore,  $\tau cl (A) = ([0,1] \cap Q) \cup A _{2}$. Again,  $f(\tau cl (A)) = \{ \frac{3}{2}, \sqrt{2} \}$;  $f ( A ) = \{ \frac{3}{2}, \sqrt{2} \}$ and $\sigma cl (f(A)) = ([1,2] \cap Q) \cup \{ \sqrt{2} \}$. So  $f ( \tau cl (A) ) \subset \sigma cl ( f (A) )$.\\
\textbf{Theorem 4.1 :} Let  a  function $f : ( X,\tau _{1},\tau _{2} ) \rightarrow ( Y, \sigma _{1}, \sigma _{2} ) $ be    pairwise  continuous  and  pairwise  open. If $A$ is  a $( i , j )$ preopen[ $( i , j )$ semi preopen] subset of $X$,  then $f(A)$ is $( i , j )$ preopen[ $( i , j )$ semi preopen] in  $Y$.\\
\textbf{Proof :} The  proof  is  similar  as  in  the  proof  of  theorem 4.1[9] and so is omitted.\\
\textbf{Theorem 4.2 :} Let  a  function $f : ( X,\tau _{1},\tau _{2} ) \rightarrow ( Y, \sigma _{1}, \sigma _{2} ) $ be    pairwise  precontinuous  and  pairwise  open. If $A$ is  a $( i , j )$ preopen[ $( i , j )$ semi preopen] subset of $Y$,  then $f^{-1}(A)$ is $( i , j )$ preopen [ $( i , j )$ semi preopen] in  $X$.\\
\textbf{Proof :} The  proof  is  similar  as  in  the  proof  of  theorem 4.2[9] and so is omitted.\\
\textbf{Theorem 4.3 :}  Let $f : ( X,\tau _{1},\tau _{2} ) \rightarrow ( Y, \sigma _{1}, \sigma _{2} ) $ be  a  function.  Then  $f$ is  pairwise  precontinuous  if  and  only  if  the  inverse  image  of  each $\sigma _{i}$ closed  set  of $Y$  is $( i , j )$ preclosed in $X$.\\
\textbf{Proof :} Let $f$ be pairwise precontinuous and let  $A$ be $\sigma _{i}$ closed  set  in $Y$.  Therefore $Y-A$ is $\sigma _{i}$ open and so $f ^{-1} (Y-A)$ is $( i , j )$ preopen  in $X$ which  implies  that $X - f ^{-1} (A)$ is $( i , j )$ preopen  in  $X$.  Hence $f ^{-1} ( A )$  is $( i , j )$ preclosed.  Conversely,  let  $A$ be  a $\sigma _{i}$ open  set  in $Y$. So $Y - A$ is $\sigma _{i}$ closed  set  in $Y$. Then  $f ^{-1} (Y-A)$ is $( i , j )$ preclosed  set  in $X$. Therefore,  $X - f ^{-1} (A)$ is $( i , j )$ preclosed in $X$. So,  $f ^{-1} ( A )$ is $( i , j )$ preopen in  $X$. Thus  $f$ is $( i , j )$ precontinuous  and hence  pairwise  precontinuous.\\
\textbf{Theorem 4.4 :} Let a function $f : ( X,\tau _{1},\tau _{2} ) \rightarrow ( Y, \sigma _{1}, \sigma _{2} ) $ be  pairwise  precontinuous. Then the following properties hold:\\
(i) for each $x \in X$ and each  $V \in \sigma _{i}$ containing $f(x)$, there exists an $ ( i , j ) $ preopen set  $U$ of $X$ containing $x$ such that $f(U) \subset V$.\\
(ii) $ f(( i , j ) pcl (A)) \subset  \sigma _{i} cl(f(A))$ for every subset $A$ of $X$.\\
(iii) $( i , j ) pcl f ^{-1} (B) \subset f ^{-1} ( \sigma _{i} cl ( B )) $ for every subset $B$ of $Y$.\\
\textbf{Proof :} (i) Let $V \in \sigma _{i}$ containing $f(x)$. So $x \in f ^{-1} (V)$. Since $f$ is  pairwise  precontinuous, $f ^{-1} (V)$ is $( i , j )$ preopen set containing  $x$. Let $U = f ^{-1} (V)$ then  $f(U) = V \subset V$.\\
(ii) Let $p \in  f(( i , j ) pcl (A))$. Clearly $A \subset  ( i , j ) pcl (A)$. Now  if $p \in f(A)$ then $p \in \sigma _{i} cl(f(A))$. So  let $p \not \in f(A)$  and $p=f ( q )$ where $q \in ( i , j ) pcl (A)$  but $q \not \in A$.  Let $U$ be  any $\sigma _{i}$ open  set  containing $p$. Then $f ^{-1} (U)$ is $( i , j )$ preopen  in $X$ containing $q$. Since $q \in ( i , j ) pcl (A)$, $A \cap f ^{-1} (U) \neq  \phi$ [by  theorem  3.6(i)]. Let $z \in A \cap f ^{-1} (U) $.  Then $f (z) \in f (A) \cap U$. But $f (z) \neq p$,  since $p \not \in f (A)$. So $U$ intersect $f (A)$ in some  point $f (z)$ other  than $p$. Then $p \in \sigma _{i} cl(f(A))$. Therefore, $ f(( i , j ) pcl (A)) \subset  \sigma _{i} cl(f(A))$.\\
(iii) Let $B$ be any subset of $Y$. Let $f ^{-1} (B) = A$. Then by (ii), $ f(( i , j ) pcl (A)) \subset  \sigma _{i} cl(f(A))$  $= \sigma _{i} cl f ( f^{-1}(B))= \sigma _{i} cl (B \cap f(X)) \subset \sigma _{i} cl(B)$. So $( i , j ) pcl  (A) \subset f ^{-1} ( \sigma _{i} cl ( B )) $ i.e., $( i , j ) pcl f ^{-1} (B) \subset f ^{-1} ( \sigma _{i} cl ( B )) $.\\
\textbf{Note 4.2 :} In  all  cases  converse  is  true  in the  case  of  bitopological  space as  seen  in  theorem 4.3[9].\\
\textbf{Theorem 4.5 :} If $f : ( X,\tau _{1},\tau _{2} ) \rightarrow ( Y, \sigma _{1}, \sigma _{2} ) $ is pairwise precontinuous and $A \in \tau _{1}\cap \tau _{2}$. Then $f : ( A,\tau _{1/A},\tau _{2/A} ) \rightarrow ( Y, \sigma _{1}, \sigma _{2} ) $ is pairwise precontinuous.\\
\textbf{Proof :} Let $U$ be a  $\sigma _{i}$ open set in $Y$. Then  $f ^{-1} (U)$ is $(\tau _{i},\tau _{j})$ preopen  in $X$. So   by theorem 3.4, $A \cap f ^{-1} (U)$ is $(\tau _{i},\tau _{j})$ preopen. Since $A \cap f ^{-1} (U) \subset A$, by theorem 3.5, $A \cap f ^{-1} (U)$ is $(\tau _{i/A},\tau _{j/A})$ preopen in $( A,\tau _{1/A},\tau _{2/A} )$. Hence $f : ( A,\tau _{1/A},\tau _{2/A} ) \rightarrow ( Y, \sigma _{1}, \sigma _{2} ) $ is  $(\tau _{i/A},\tau _{j/A})$  precontinuous and hence pairwise precontinuous.\\
\textbf{Theorem 4.6 :} Let $f : ( X,\tau _{1},\tau _{2} ) \rightarrow ( X^{*}, \tau _{1} ^{*}, \tau _{2}^{*} ) $ be  a $( i , j )$ precontinuous  mapping  such  that  the  following  condition 'C' is satisfied.\\
\indent C :  $f ( \tau _{j} cl f ^{-1} ( U ^{*}))  = U^{*}$ for  every $\tau _{i} ^{*}$ open  set $U^{*} $. \\
Then if $\{ x _{\alpha} : \alpha \in D \}$ is  a  net  converging  to $x$ in $( X, \tau _{i} )$ then $\{ f ( x _{\alpha} ) : \alpha \in D \}$ converges  to $f ( x )$ in  $( X ^{*}, \tau _{i} ^{*} )$.\\
\textbf{Proof :} Let $U ^{*}$ be a $\tau _{i} ^{*}$ open set  in $X ^{*}$  containing $f (x)$. So $x \in f ^{-1} ( U ^{*} )$. Since $f$  is $( i , j )$ precontinuous, $ f ^{-1} ( U ^{*} )$  is $( i , j )$ preopen. So  there  exists  a $\tau _{i} $ open  set $U$ such  that $x \in f ^{-1} ( U ^{*} ) \subset U \subset \tau _{j} cl f ^{-1} ( U ^{*}) $.  Since $\{ x _{\alpha} : \alpha \in D \}$ converges  to $x$ in $( X, \tau _{i} )$ there  exists $\alpha _{o} \in D$ such  that $x _{\alpha} \in U$, for  all $\alpha \geq \alpha _{o}$. So  $f( x _{\alpha}) \in f (U) \subset f ( \tau _{j} cl f ^{-1} ( U ^{*}))  = U^{*}$  [ by the condition C ]. Hence  $\{ f ( x _{\alpha} ) : \alpha \in D \}$  is  eventually  in  every  open  set $U ^{*}$ in $( X ^{*}, \tau _{i} ^{*} )$. So $\{ f ( x _{\alpha} ) : \alpha \in D \}$ converges  to $f (x)$ in  $( X ^{*}, \tau _{i} ^{*} )$.\\

\noindent\textbf{5. Pairwise semi precontinuity ( sp-continuity ):} \\
\textbf{Definition 5.1 ( cf. [9] ) :} A function $f : ( X,\tau _{1},\tau _{2} ) \rightarrow ( Y, \sigma _{1}, \sigma _{2} ) $ is  said  to  be pairwise semi precontinuous or  pairwise  sp-continuous ( resp. pairwise semi continuous ) if the inverse image of each  $\sigma _{i}$ open set of $Y$ is $( i , j )$ semi preopen ( resp. $( i , j )$ semi open ) in $X$ where $i \neq j, i,j=1,2.$\\
\textbf{Remark 5.1 :} As in the case of bitopological space [9], it can be easily checked that pairwise continuity implies pairwise semi continuity and pairwise precontinuity. Also pairwise semi continuity implies pairwise sp-continuity and pairwise precontinuity implies pairwise sp-continuity. However the reverse implications are not true in a bispace even in case of a bitopological space also which is  seen in [9] \\
\textbf{Theorem 5.1 :} Let $f : ( X,\tau _{1},\tau _{2} ) \rightarrow ( Y, \sigma _{1}, \sigma _{2} ) $ be a function. Then  $f$ is pairwise sp-continuous if and only if  the inverse image of each $\sigma _{i}$ closed set of $Y$ is $( i , j )$ semi preclosed in  $X$.\\
\textbf{Proof :} Proof is similar to the proof of theorem 4.3.\\
\textbf{Theorem 5.2 :} Let a function $f : ( X,\tau _{1},\tau _{2} ) \rightarrow ( Y, \sigma _{1}, \sigma _{2} ) $ be pairwise sp-continuous.  Then the following properties hold.\\
(i) for each $x \in X$ and each  $V \in \sigma _{i}$ containing $f(x)$, there exists an $( i , j )$ semi preopen set  $U$ of $X$ containing $x$ such that $f(U) \subset V$.\\
(ii) $ f(( i , j ) spcl (A)) \subset  \sigma _{i} cl(f(A))$ for every subset $A$ of $X$.\\
(iii) $( i , j ) spcl f ^{-1} (B) \subset f ^{-1} ( \sigma _{i} cl ( B )) $ for every subset $B$ of $Y$.\\
\textbf{Proof :}  Proof is similar to the proof of theorem 4.4.\\
\textbf{Theorem 5.3 :} If $f : ( X,\tau _{1},\tau _{2} ) \rightarrow ( Y, \sigma _{1}, \sigma _{2} ) $ is pairwise sp-continuous and $A \in \tau _{1}\cap \tau _{2}$. Then $f : ( A,\tau _{1/A},\tau _{2/A} ) \rightarrow ( Y, \sigma _{1}, \sigma _{2} ) $ is pairwise sp-continuous.\\
\textbf{Proof :} Proof is similar to the proof of theorem 4.5.\\

\begin{center}
{\bf\Large References}\\
\end{center}
1) Alexandroff.A.D, \emph{Additive set functions in abstract spaces}, (a) Mat.Sb.(N.S),8:50(1940),\\ \indent 307-348.(English,Russian Summary). (b) ibid,
    9:51(1941), 563-628, (English,Russian Summary).\\
(2) Banerjee A.K. and Saha P.K., \emph{Bispace Group,} Int.J.of Math. Sci.and Engg. Appl.(IJMESEA);  Vol.5 No V(2011) pp. 41-47.\\
(3) Banerjee A.K. and Saha P.K., \emph{Semi open sets in bispaces}, CUBO A Mathematical Journal; Vol. 17, No. 01,(2015),99-106.\\
(4) Banerjee A. K. and Saha P. K., \emph{Quasi open sets in bispaces}, General Math. Notes, Vol. 30, No. 2, (2015), pp. 1-9.\\
(5) Bose S.,  \emph{Semi open sets, Semi continuity and Semi open mapping in bitopological spaces}, Bull. Cal. Math. Soc. 73,(1981), 237-246 . \\
(6) Jelic M., \emph{A decomposition of pairwise continuity}, J. Inst. Math. Comput. Sci. Math. Soc., 3, (1990), 25-29 .\\
(7) Kar A. and Bhattacharyya P., \emph{Bitopological preopen sets, Precontinuity and Preopen Mappings}, Indian J. of Math., 34 ,(1992), 295-309.\\
(8) Kelley J.C.,\emph{ Bitopological spaces}, Proc.London.Math.Soc., 13(1963), 71-89.\\
(9) Khedr F. H. and Areefi S. M., \emph{Precontinuity and semi precontinuity in bitopological spaces}, Indian J. of Pure appl. Math., 23(9),(1992), 625-633.\\
(10) Lahiri B.K. and Das P., \emph{Semi Open set in a space}, Sains malaysiana 24(4), (1995), 1-11.\\
(11) Lahiri B.K. and Das P., \emph{Certain bitopological concept in a Bispace}, Soochow J.of Math.Vol.27.No.2,(2001),
 pp.175-185.\\
(12) Levine N., \emph{Semi Open sets and Semi Continuity in topological spaces}, Amer.Math.Monthly
\indent 70(1963), 36-41.\\
(13) Mashhour A. S.,  Monsef M. E. Abd. El. and Deeb S. N. E., \emph{On precontinuous and weak precontinuous mappings}, Proc. Math. Phys. Soc. Egypt 53,(1982), 47-53.\\
(14) Pervin W.J., \emph{Connectedness in Bitopological spaces}, Ind.Math.29,(1967), 369-372.\\
(15) Reilly I. L., \emph{On bitopological separation properties}, Nanta Math., 5 (1972), 14-25.\\

\noindent $^{1}$ Department of Mathematics, The University of Burdwan \\
Burdwan-713104, W.B., India\\
Email: akbanerjee@math.buruniv.ac.in\\
$^{2}$ Department of Mathematics, Behala College \\
Kolkata-700060, W.B., India\\
Email: pratapsaha2@gmail.com\\

\end{document}